\documentclass[12pt, reqno]{amsart}
\usepackage{amsmath, amsthm, amscd, amsfonts, amssymb, graphicx, color}
\usepackage[bookmarksnumbered, colorlinks, plainpages]{hyperref}

\usepackage[all]{xy}
\newtheorem{theorem}{Theorem}[section]

\newtheorem{proposition}[theorem]{Proposition}
\newtheorem{corollary}[theorem]{Corollary}
\theoremstyle{definition}

\theoremstyle{remark}

\numberwithin{equation}{section}

\begin{document}

	\setcounter{page}{1}
	\title[Reciprocal Dunford--Pettis sets and V$^*$-sets in Banach lattices]{Reciprocal Dunford--Pettis sets and V$^*$-sets in Banach lattices}
	
\author[J.X. Chen]
       {Jin Xi Chen}
       \address{School of Mathematics, Southwest Jiaotong University, Chengdu 610031, China}
       \email{jinxichen@swjtu.edu.cn}

	 \author[X. Li]
	  {Xi Li}
	 \address { School of Mathematics, Southwest Jiaotong University, Chengdu 610031, China}
	
\email{lixi@my.swjtu.edu.cn}

	

	\begin{abstract}
   In this short note, we show that one cannot differentiate between reciprocal Dunford--Pettis sets and V$^*$-sets in a Banach lattice. That is, for a bounded subset $K$ of a Banach lattice $E$, $K$ is a V$^*$-set if and only if $K$ is a reciprocal Dunford--Pettis set, or equivalently, if and only if every disjoint sequence in  the solid hull of $K$ is weakly null (i.e. $K$ is disjointly weakly compact).

	\end{abstract}
\maketitle\baselineskip 4.90mm
Throughout this note, $X, Y$ will denote real Banach spaces, and $E, F$ will denote real Banach lattices. $B_X$ is the closed unit ball of $X$.  The solid hull of a subset $K$ of $E$ is denoted by $Sol(K):=\lbrace y\in E:|y| \leq |x|$ \ for \ some\ $x\in K\rbrace$.

\par Following A. Pe{\l}czy\'{n}ski  \cite{Pel}, a bounded subset $K$ of  $X$ (resp. $K$ of $X^\prime$) is called a \textit{V$^{\,*}$-set} (resp. \textit{V-set}) if
     $$\sup_{x\in K}|x_{n}^{\,\prime}(x)|\to0 \,\,\,\,\,\,\textmd{(}\textmd{resp}. \sup_{x^{\,\prime}\in K}|x^{\,\prime}(x_n)|\to0 \textmd{)}\,\,\, \,\, (n\to\infty)$$for each weakly unconditionally Cauchy series $\sum_{n} x^{\,\prime}_{n}$ in $X^\prime$  (resp. $\sum_{n} x_{n}$ in $X$).


 \par Recall that a Banach space $X$ is said to have the \textit{Dunford--Pettis property} if every weakly compact operator from $X$ into an arbitrary Banach space is Dunford-Pettis, and $X$ is said to have the \textit{reciprocal Dunford--Pettis property} if every Dunford-Pettis operator from  $X$ into an arbitrary Banach space is weakly compact. This implies that  $X$ has the reciprocal Dunford--Pettis property if and only if $TB_X$ is relatively weakly compact for every Dunford--Pettis operator $T$ with domain $X$. Inspired by this observation, Bator, Lewis and  Ochoa \cite{Bator} introduced a class of sets that they called reciprocal Dunford--Pettis sets. A bounded subset $K$ of $X$ is called a \textit{reciprocal Dunford--Pettis set} (or \textit{RDP}-set) if $T(K)$ is relatively weakly compact for every Banach space $Y$ and every Dunford--Pettis operator $T:X\to Y$. Bator et. al \cite[Theorem 4.2]{Bator} established the following containment relationships that exist among the classes of sets mentioned above:
$$\textmd{WPC}\subseteq \textmd{RDP}\subseteq \textmd{V}^*$$Here,  WPC stands for the family of weakly precompact sets (i.e.  weakly conditionally compact sets). Furthermore, each containment is proper. However,  Bator et. al \cite[p.10]{Bator} and Bombal \cite{Bom} realized that one cannot use $C(\Omega)$ or $C(\Omega)^\prime$ to differentiate between RDP-sets and V$^*$-sets since the RDP-sets and the V$^*$-sets coincide in these two types of spaces. Note that both $C(\Omega)$ and  $C(\Omega)^\prime$ are  Banach lattices.
\par The purpose of this note is to show that one cannot expect to differentiate between reciprocal Dunford--Pettis sets and V$^*$-sets in the Banach lattice context. That is, for a bounded subset $K$ of a Banach lattice $E$, $K$ is a V$^*$-set if and only if $K$ is a reciprocal Dunford--Pettis set.

\par  We  employ  disjointly weakly compact sets to bridge the gap between reciprocal Dunford--Pettis sets and V$^*$-sets.  Following W. Wnuk \cite{WOrder}, a bounded subset $K$ of a Banach lattice $E$ is called a \textit{disjointly weakly compact set} (or \textit{DWC-set}) whenever $x_{n}\xrightarrow{w}0$ for every disjoint sequence  $(x_n)\subset Sol(K)$.  Recently, Xiang, Chen and Li \cite{XCL, disjointly weakly compact set} extensively studied  the properties  of disjointly weakly compact sets and related operators. Also, a  bounded linear operator $T: E\rightarrow X$ from a Banach lattice into a Banach space  is called an \textit{almost Dunford--Pettis operator} if $\|Tx_{n}\| \rightarrow 0$ for each disjoint weakly null sequence $(x_n)$ in $E$ \cite{Sanchez}. Clearly, every Dunford--Pettis operator on a Banach lattice is almost Dunford--Pettis. 

\par We are now in a position to show that  RDP-sets, V$^{\,*}$-sets and disjointly weakly compact sets always coincide in a Banach lattice. Although the property (V$^*$) of Banach lattices were disscussed in  some papers (e.g., \cite{Niculescu, Saab}), or mentioned in  monographs \cite[Corollary 5.3.5]{Meyer} and \cite[Theorem 7.3]{WOrder}, the result on the coincidence of the V$^*$-sets and the disjointly weakly compact sets has never appeared in literature as far as the authors know.
\begin{proposition} \label{relationship} For a bounded subset $K$ of a Banach lattice $E$ the following statements are equivalent:
\begin{enumerate}
  \item $K$ is an RDP-set.
  \item $K$ is a V$^{\,*}$-set.
  \item $K$ is a disjointly weakly compact set.
  \item $T(K)$ is relatively weakly compact for every Banach space $X$ and every almost Dunford-Pettis operator $T:E\to X$.
\end{enumerate}
 \end{proposition}

 \begin{proof}
 \par$(1)\Rightarrow(2)$ This follows from Theorem 4.2 of \cite{Bator}.
 \par$(2)\Rightarrow(3)$ Assume that $K$ is a  V$^{*}$-subset of $E$. To prove that $K$ is disjointly weakly compact, we have to show that every order bounded disjoint sequence from $E^{\,\prime}$ converges uniformly to zero on $K$ (see \cite[Theorem 2.3]{disjointly weakly compact set}). To this end, let  $(x^{\,\prime}_{n})\subset E^{\,\prime}$ be an order bounded disjoint sequence. Then the series $\sum_{n} x^{\,\prime}_{n}$ is weakly unconditionally Cauchy \cite[p.192]{Positive}. Therefore, we have $\sup_{x\in K}|x_{n}^{\,\prime}(x)|\to 0$.
 \par$(3)\Rightarrow(4)$   Let $T:E\to X$ be an almost Dunford-Pettis operator from $E$ into an arbitrary Banach space $X$. By definition, we can assume without loss of generality that $K$ is solid. For every disjoint sequence $(x_n)\subset K$, we have $x_{n}\xrightarrow{w}0$, and hence $\|Tx_n\|\to 0$ since $T$ is almost Dunford-Pettis. Then in view of Theorem 4.36 of \cite{Positive}, for each $\varepsilon>0$  there exists some $u\in{E^{+}}$ lying in the ideal generated by $K$ such that
	$$\|T[(|x|-u)^{+}]\|<\frac{\varepsilon}{2}$$holds for all $x\in K$. Therefore, from the identities $x^{+}=x^{+}\wedge u+(x^{+}-u)^+$ and $x^{-}=x^{-}\wedge u+(x^{-}-u)^+$ we can see that $x\in [-u,u]+(x^{+}-u)^{+}-(x^{-}-u)^+$. Thus, we have
$$T(K)\subset T[-u,u]+\varepsilon B_{X}$$
Since $T$ is almost Dunford-Pettis,  $T$ is an order weakly compact operator, and hence  $T[-u,u]$ is relatively weakly compact. It follows that  $T(K)$ is a relatively weakly compact subset of $X$.
\par$(4)\Rightarrow(1)$ Obvious.
 \end{proof}

 Combined Proposition \ref{relationship} with Theorem 2.4 of \cite{XCL}, we have indeed established the following containment relationships in the Banach lattice context:
$$\textmd{WPC}\subseteq \textmd{DWC}=\textmd{RDP}= \textmd{V}^*$$

 As an immediate consequence of Proposition \ref{relationship} the following result shows that a bounded linear operator between two Banach lattices always preserves the disjointly weak compactness of a set.

\begin{corollary}
Let $T:E\to F$ be a continuous linear operator between two Banach lattices. Then $T(K)$ is likewise disjointly weakly compact for every disjointly weakly compact subset $K$ of $E$.
\end{corollary}

 Recall that a Banach space $X$ is said to have the \textit{property}
 (\textit{V$^{\,*}$}) of Pe{\l}czy\'{n}ski if every V$^*$-set in $X$ is relatively weakly compact, and $X$ is said to have the \textit{property weak} (\textit{V$^{\,*}$}) if every V$^*$-set in $X$ is  weakly precompact \cite[p.111]{Bom2}. In view of our alternative characterization of V$^*$-sets in Banach lattices, some known results related to V$^*$-sets can be reformulated and re-proved in a simpler way with no function theory involved. We take some examples:
 \begin{enumerate}
   \item A Banach lattice $E$ has the property (V$^{\,*}$) if and only if $E$ is weakly sequentially complete (i.e. $E$ is a $KB$--space) \cite[Proposition 2.6] {disjointly weakly compact set}.
   \item A Banach lattice $E$ has the property weak (V$^{\,*}$) if and only if every order interval in $E$ is weakly precompact \cite[Theorem 2.4]{XCL}.
 \end{enumerate}

 \begin{corollary} \cite[Theorem 2.2]{Bom2} Let $E$ be a Banach lattice with order continuous norm. Then every closed subspace of $E$ has the
 property weak \textmd{(V$^{\,*}$)}.
 \end{corollary}

 \begin{proof}
 Let $X$ be a closed subspace of $E$ and let $K$ be a V$^*$-subset of $X$. Then $K$ is also a V$^*$-subset of $E$. By Proposition \ref{relationship} $K$ is disjointly weakly compact in $E$. Since $E$ has order continuous norm (and hence every order interval in $E$ is weakly compact), it follows from \cite[Theorem 2.4]{XCL} that $K$ is weakly precompact in $E$, and hence $K$ is weakly precompact in $X$. This implies that $X$ has the
 property weak \textmd{(V$^{\,*}$)}.
 \end{proof}

Naturally, one may ask whether V-sets in the dual of a Banach lattice also possess some type of disjoint weak$^*$ compactness property. For our convenience, let us define a bounded subset $A$ of $E^\prime$ to be a \textit{disjointly weak$^{\,*}$  compact} if $x_{n}^{\,\prime}\xrightarrow{w^*}0$ for every disjoint sequence  $(x_{n}^{\,\prime})\subset Sol(A)$. It is well known that a Banach lattice $E$ has order continuous norm if and only if $B_{E^{\prime}}$ is disjointly weak$^{*}$  compact (see, e.g., \cite[Corollary 2.4.3]{Meyer}). Now, let $A$ be a V-subset of $E^\prime$. Then, for every order bounded disjoint sequence $(x_n)$ of $E$, we have  $\sup_{x^{\,\prime}\in A}|x^{\,\prime}(x_n)|\to 0$  since the series $\sum_{n}x_n$ is weakly unconditionally Cauchy. From a result due to P. G. Dodds and D. H. Fremlin \cite{1979Compact},   it follows that $A$ is disjointly weak$^{*}$  compact (see, e.g., \cite[Theorem 5.63]{Positive}). However, a disjointly weak$^*$ compact set is not necessarily a V-set. Indeed,  $B_{\ell^1}$ is an example of a disjointly weak$^{*}$  compact set which is not a V-set.

\end{document}